\documentstyle[]{article}

\begin{document}

\title{Fractional Fokker--Planck Equation for  Nonlinear Stochastic
Differential Equations Driven by Non-Gaussian Levy Stable Noises}

\author{D. Schertzer,
M. Larchev\^eque \\
Laboratoire de Mod\'elisation en M\'ecanique, Tour 66, Boite 162
\\ Universit\'e
Pierre et Marie Curie,\\ 4 Place Jussieu, F-75252 Paris Cedex 05, France.\\
\and
J. Duan \\ Department of Mathematical Sciences \\
Clemson University \\
 Clemson, SC 29634, USA.\\
\and
V.V. Yanovsky \\
Turbulence Research, Institute for Single Crystals\\
 National Acad. Sci. Ukraine, \\ Lenin ave. 60, Kharkov 310001, Ukraine\\
\and
S. Lovejoy \\ Physics Department, McGill University \\
3600 University Street
Montreal, H3A 2T8, Quebec, Canada.}

\date{Oct 28 1999}

\maketitle


\def\u{\underline}
\def\a{\alpha}\def\b{\beta}\def\g{\gamma}\def\d{\delta}\def\D{\Delta}
\def\e{\epsilon}\def\k{\kappa}\def\l{\lambda}\def\L{\Lambda}
\def\s{\sigma}\def\S{\Sigma}\def\Th{\Theta}\def\th{\theta}
\def\om{\omega}\def\Om{\Omega}\def\G{\Gamma}\def\v{\varepsilon}

\newtheorem{theorem}{Theorem}
\newtheorem{definition}{Definition}
\newtheorem{proposition}{Proposition}
\newtheorem{cor}{Corollary}
\newtheorem{lem}{Lemma}

\abstract{{\bf J. Math. Phys., 42(2001), 200-212.}
The
Fokker-Planck equation has been very useful for studying dynamic
behavior of stochastic differential equations driven by Gaussian
noises. In this paper, we derive a Fractional Fokker--Planck
equation for the probability distribution of particles whose
motion is governed by a {\em nonlinear} Langevin-type equation,
which  is driven by a non-Gaussian Levy-stable noise. We obtain in
fact a more general result for Markovian processes generated by
stochastic differential equations. }


Correspondence should be addressed to D. Schertzer (fax: +33 1 44275259,
e-mail: schertze@ccr.jussieu.fr).
\newpage

\section{Introduction and statement of the problem} \label{introduction}
The Fokker-Planck equation is one of the most celebrated equations in
Physics, since
it has been very useful for studying
dynamic behavior of stochastic differential equations
driven by Gaussian noises.
There has recently been a mushrooming interest \cite{ Fogedby, Chechkin,
Compte, Marsan, Chaves, Yan}
in the fact that the probability
density of a (linear) Levy motion satisfies a generalized Fokker-Planck
equation involving
fractional orders of differentiation.
This essentially corresponds to a re-interpretation of the characteristic
function of a Levy
motion and significant applications seem to require its nonlinear
generalization.

We therefore consider the following nonlinear Langevin--like equation for a
stochastic (real) quantity $X(t):$

\begin{equation}  \label{eq_langevin}
d X(t) = m(X(t),t)~ dt + \sigma(X(t),t)~ dL
\end{equation}

\noindent where the driving source is a Levy stable motion $L$, i.e.
a motion (e.g. \cite{Fristedt}) whose increments $\delta L$ are
stationary and independent
for any time lag $\delta t$ and correspond to independent, identically
distributed Levy stable variables
\cite{Levy,Khintchine, Gnedenko, Feller, Zolotarev}. Let
us recall that a  Levy stable motion is defined, as are its increments, by
four parameters:
its Levy stability index $\alpha$ ($0 <\alpha \leq 2$) , its skewness
$\beta$ ($-1 \leq \beta \leq 1$),
its center $\gamma \delta t$ and its scale parameter $ D~ \delta t$ ($D
\geq 0$).
Brownian motion corresponds to the limit case $\alpha = 2$, which also
implies $\beta = 0 $, and to the 'normal' diffusion law. The variance
$Var[X(t) -X(t_0)]$
of the distance traveled  by a brownian particle,
is twice its scale parameter and therefore yields the classical Einstein
relation:
$ Var[X(t) -X(t_0)] = 2 D (t-t_0)$.

The linear case, which is the unique case studied until now, corresponds to:

\begin{equation}  \label{eq_linear_case}
m(X(t),t) \equiv m = Const., ~~ \sigma(X(t),t) \equiv \sigma = Const.
\end{equation}

\noindent $X(t) -X(t_0)$ is also a Levy motion which has the same Levy
stability
index $\alpha $,
but with a possible different center or trend
(when $ m \neq 0$) and scale or amplitude (when $\sigma \neq 1$).

In the nonlinear cases,
$\sigma (X(t),t)\geq 0$ and $m(X(t),t)$ are nonlinear
functions of $X(t)$ and $t$, which satisfy certain regularity constraints to
be discussed later.
We claim that:

\begin{proposition}\label{prop.1}
\begin{quote} The transition probability density:

\begin{equation} \label{eq.trans_proba}
\forall t\geq t_{0}:~~p(x,t|x_{0},t_{0})=Pr(X(t)=x|X(t_{0})=x_{0})
\end{equation}

\noindent corresponding to the nonlinear stochastic differential equation
(Eq.\ref{eq_langevin}), with $\alpha \neq 1 ~or~ \beta = 0$,
 is the solution of the
following Fractional Fokker-Planck equation: .

\begin{eqnarray}  \label{eq_fokker_planck}
{\frac{\partial }{\partial t}} p(x,t|x_0,t_0) = - {\frac{{\partial} }{{%
\partial x}}}(\gamma \sigma (x,t) + m(x,t)) p(x,t|x_0,t_0)   \nonumber \\
 -  D [(-\Delta )^{\alpha /2} +
\beta\omega (\alpha ) {\frac{\partial }{\partial x}} (-\Delta )^{(\alpha
-1)/2}] \sigma (x,t)^{\alpha} p(x,t|x_0,t_0)
\end{eqnarray}

\end{quote}
\end{proposition}

\noindent where $\omega(\alpha)$ is defined by:

\begin{equation}  \label{eq_omega}
\alpha \neq1:~~ \omega (\alpha ) = tan {\frac{\pi \alpha }{2}}
\end{equation}

\noindent  and where the fractional powers of the Laplacian $\Delta$ will
be discussed
in Sect.\ref{sec.fractional}. Proposition \ref{prop.1} and
Eq.\ref{eq_fokker_planck}
will be established for scalar processes
(i.e. $\Delta \equiv  {\frac{\partial^2 }{\partial x^2}}$)
and its extension to vector processes will be discussed and presented in
Sect.\ref{sec.high.dim}.

This Fractional Fokker-Planck equation
will be established with the help of the much more general proposition:

\begin{proposition}\label{prop.2}
\begin{quote}
The inverse Fourier transform of the second characteristic
function or cumulant generating function of the increments of a Markov
process $X(t)$ generates by convolution the Fokker-Planck equation of
evolution of its
transition probability $p(x,t|x_{0},t_{0})$.
\end{quote}
\end{proposition}

We will demonstrate this proposition in a straightforward, yet rigorous way.
More precisely, we will establish the
following:

\begin{equation} \label{eq.FP.generation}
{\frac{\partial p}{\partial t}}(x,t|x_{0},t_{0}) =
\int {dy {\frac{\partial \widetilde{K}}{\partial t}}(x-y|y,t)
p(y,t|x_{0},t_{0})}
\end{equation}
\noindent where $\widetilde{K}$ is the inverse Fourier transform of the
cumulant generating
function of the increments. Its arguments will become explicit in
Sect.\ref{cumulants.sect}.

This not only holds for
processes with stationary and independent increments,
as in the linear case (Eq.\ref{eq_linear_case})
but for any Markov process, including those defined by
the non-linear Langevin-like equation (Eq.\ref{eq_langevin}
with $m \ne Const.$, $\sigma \ne Const.$).
As a consequence of Eq.\ref{eq.FP.generation}, we will demonstrate the
following:

\begin{proposition}\label{prop.3}
\begin{quote}
The Kramers-Moyal coefficients $A_{n}$ of the
Fokker-Planck equation of a Markov process $X(t)$:

\begin{equation}\label{eq_NFP_markov.3}
{\frac{\partial p}{\partial t}}(x,t|x_{0},t_{0})=\sum_{n \in J}
  {\frac{\partial ^{n}}{\partial x^{n}}}
[A_{n}(x,t) p(x,t|x_{0},t_{0})]
\end{equation}

\noindent are directly related to the cumulants $C_{n}$ of the increments:

\begin{equation}
A_{n}(x,t) = {\frac{(-1)^n} {n!}} C_{n}(x,t)
\end{equation}

\end{quote}

\end{proposition}

\noindent where the set $J$ of the indices $n$ is $\{1, 2\}$ in
the most classical case (e.g. \cite{Gardiner}, which is a particular case of
$J \subseteq {\bf N }$ which corresponds to an analytic expansion of cumulants.

We will demonstrate this property (Prop.\ref{prop.3}),
which at best is only mentioned in few standard text books on the
(classical) Fokker-Plank equation, as well as its generalization for
non analytic cumulant expansions, i.e. there are non integers indices ${n
\in J}$.
This latter property, discussed in Sect.\ref{sec.fractional},
will be exploited in Sect.\ref{Levy case} in order to derive
Prop.\ref{prop.1}
with $J= \{1, \alpha \}, 0< \alpha \le 2$.

\section{The cumulant generating function of the increments}
\label{cumulants.sect}

The first and second (conditional) characteristic functions are
respectively the moment generating function $Z_{X}(k,t-t_{0}|x_{0},t_{0})$ and
the cumulant generating function $K_{X}(k,t-t_{0}|x_{0},t_{0})$, associated
with the transition probability $p(x,t|x_{0},t_{0})$ of a process $X(t)$.
They are defined
by the Fourier transform of the latter, with $k$ being the conjugate
variable of $x-x_{0}$ :

\begin{eqnarray}  \label{eq_fc_Z}
F[p(x,t|x_{0},t_{0})] &\equiv & Z_{X}(k,t-t_{0}|x_{0},t_{0})  \\
&\equiv & exp(K_{X}(k,t-t_{0}|x_{0},t_{0})) \\
&\equiv & E [exp(ik(X(t)-X_{0}))| X(t_{0})] = X_{0}]
\end{eqnarray}

\noindent where  $E[ \cdot | \cdot ]$ denote the conditional mathematical
expectation,
$F$ and $F^{-1}$  respectively the Fourier--transform
and its inverse:

\begin{eqnarray}
F[f] &=&\hat{f}(k)=\int_{-\infty }^{\infty }{dx}~~exp(ikx)f(x)~~~
\label{eq.fourier} \\
~F^{-1}[\hat{f}] &=&f(x)=\int_{-\infty }^{\infty }{\frac{dk}{2\pi }}%
~~exp(-ikx)\hat{f}(k)
\end{eqnarray}

The corresponding quantities for increments $\delta X(\delta t)=X(t +
\delta t)-X(t)$,
corresponding to a given time lag $\delta t>0$, are defined in a
similar way:

\begin{eqnarray}
F[p(x +\delta x , t+\delta t|x,t)] &=&\delta Z_{X}(k,\delta t|x,t)
\label{eq_incr} \\
&\equiv & exp(\delta K_{X}(k,\delta t|x,t)) \\
& = & E[ exp(ik (X(t + \delta t)-X))|  X(t)= X]
\end{eqnarray}

\noindent where $k$ is the conjugate variable of $\delta x$.
The cumulants of the increments $C_{n}$
are the coefficients of the
Taylor expansion of ${\delta K_{X}}$:

\begin{equation} \label{eq.cumulant_gen}
{\delta K_{X} (k,\delta t|x,t)=\delta t}
\sum_{n \in J}{\frac {(ik)^{n}} {n!}}C_{n}(x,t)+o({\delta t)}
\end{equation}

As already mentioned, the classical case corresponds to
an analytic expansion of ${\delta K_{X}}$, i.e. $J \subseteq {\bf N}$, whereas
we will be interested by a non-analytic expansion $J= \{1, \alpha \}$.

\section{Processes with stationary and independent increments}
\label{sect.sii_porcesses}
Let us first consider the simple sub-case of a process with
stationary and independent increments. It corresponds to
$ C_{n}(x,t) \equiv  C_{n} = Const.$ in Eqs.\ref{eq_NFP_markov.3},
\ref{eq.cumulant_gen} and as already discussed in Sect. \ref{introduction},
it includes the linear case (Eq.\ref{eq_linear_case}) of the Langevin--like
equation (Eq.\ref{eq_langevin}).

However we believe that the following derivation is not
only somewhat pedagological on the role of the characteristic functions
for the nonlinear case, but terser than
derivations previously presented for the linear case.

The stationarity of the increments implies that the transition probability
depends only on the time and
space lags, i.e.:

\begin{equation} \label{eq.pdf_ind_incr}
p(x,t|x_{0},t_{0})=p(x-x_{0},t-t_{0})
\end{equation}

\noindent and similarly, the characteristic
functions of the increments are no longer conditioned, for instance:

\begin{eqnarray}
Z_{X}(k,t-t_{0}|x_{0},t_{0}) &\equiv & Z_{X}(k,t-t_{0})  \\
K_{X}(k,t-t_{0}|x_{0},t_{0}) &\equiv & K_{X}(k,t-t_{0})
\end{eqnarray}

On the other hand, the independence of the increments
implies that the transition probabilities satisfy a convolution (over any
possible intermediate
position $y$) for any  given time lag $\delta t$:

\begin{equation} \label{eq.convolution}
\forall  \delta t >0 : p(x-x_{0},t+ \delta t-t_{0}) = \int dy ~~
p(x-y,\delta t)p(y - x_{0},t-t_{0})
\end{equation}

\noindent and the corresponding characteristic functions merely
factor (resp. add). Therefore, we have:

\begin{equation} \label{eq.fc_ind_incr}
Z_{X}(k,t+\delta t-t_{0})-Z_{X}(k,t-t_{0})=
Z_{X}(k,t-t_{0})({\delta Z_{X}(k,\delta t)-1)}
\end{equation}

This in turn leads to:

\begin{equation} \label{eq.fourier_chap_kolmo}
Z_{X}(k,t+\delta t-t_{0})-Z_{X}(k,t-t_{0})=
Z_{X}(k,t-t_{0}) {\delta K_{X}(k,\delta t)} + o(\delta t)
\end{equation}

Its inverse Fourier transform yields:

\begin{equation} \label{eq_NFP_ind.inc.1}
p(x,t+\delta t|x_{0},t_{0})-p(x,t|x_{0},t_{0}) =
\int {d y} F^{-1}[\delta K_{X}(k,\delta t)]
p(y-x_{0},t-t_{0}) + o(\delta t)
\end{equation}

This demonstrates (in the limit $\delta t\to 0$)
Prop.\ref{prop.2} and Eq.\ref{eq.FP.generation},
as well as Prop.\ref{prop.3}, since Eq.\ref{eq_NFP_ind.inc.1} corresponds,
with the help of Eq.\ref{eq.cumulant_gen}, to:

\begin{equation} \label{eq_NFP_ind.inc.2}
p(x,t+\delta t|x_{0},t_{0})-p(x,t|x_{0},t_{0}) =\delta t\sum_{n}
[C_{n}{\frac {(-1)^{n}} {n!}}
\int dy\delta _{x-y}^{(n)} p(y,t|x_{0},t_{0})] + o(\delta t)
\end{equation}

\noindent  where $\delta _{x}^{n}$ denotes
the $n^{th}$ derivative of the Dirac function. Therefore, we obtain:

\begin{equation} \label{eq_NFP_ind.inc.3}
{\frac {\partial }{\partial t}} p(x,t|x_0,t_0) = \sum_{n \in J}
A_{n}{\frac{\partial ^{n}}{\partial x^{n}}} p(x,t|x_{0},t_{0})
\end{equation}

\noindent which corresponds to the linear case of Eq.\ref{eq_NFP_markov.3}.

\section{More general Markov processes} \label{sect.Markov}

In the case of a  Markov process which does not have stationary and
independent increments,
there is no longer a simple convolution equation (Eq. \ref{eq.convolution}) of
the transition probabilities, nor
a simple factorization of characteristic functions
(Eq.\ref{eq.fc_ind_incr}).
However, the former satisfies a generalized convolution equation which
corresponds to the
Chapman-Kolmorogorov identity \cite{Feller} valid
for any Markov process $X(t)$:

\begin{equation} \label{eq.chap_kolmo}
\forall \delta t > 0:  p(x,t+\delta t|x_{0},t_{0}))=\int {d y
~p(x,t+\delta t|y,t)p(y,t|x_{0},t_{0})}
\end{equation}

\noindent which indeed reduces to a mere convolution  (Eq.
\ref{eq.convolution})
in the case of processes with
stationary and independent increments. This identity can be written under
the equivalent form:

\begin{equation}
p(x,t+\delta t|x_0,t_0) = \int{d y \int{{\frac{dk }{2\pi }} e^{-i k
y + \delta K_{X}(k,\delta t |y,t )}} p(y,t|x_0,t_0)}
\end{equation}

Noting that we have:

\begin{equation}
p(x,t|x_{0},t_{0}) =\int {dy~p(y,t|x_{0},t_{0})}%
\int {{\frac{dk}{2\pi }}e^{-ik y}}
\end{equation}

\noindent we obtain:

\begin{equation} \label{eq_NFP_markov.1}
p(x,t+\delta t|x_{0},t_{0})-p(x,t|x_{0},t_{0}) = \delta t \int {d y}
F^{-1}[\delta K_{X}(k,\delta t|y,t)]
p(y,t|x_{0},t_{0}) + o(\delta t)
\end{equation}

In the limit $\delta t\to 0$, this corresponds to
Prop. \ref{prop.2} and Eq. \ref{eq.FP.generation}.
When $J \subseteq {\bf N}$, it yields with the help of
Eq.\ref{eq.cumulant_gen}:

\begin{equation} \label{eq_NFP_markov.2}
\delta p(x,t|x_{0},t_{0})=\delta t\sum_{n \in N} \int {d y} \delta
_{x-y}^{(n)}
[{\frac {(-1)^{n}} {n!}}C_{n}(y,t)p(y,t|x_{0},t_{0})]+o(\delta t)
\end{equation}

The limit $\delta t \to 0$
corresponds to Eq.\ref{eq_NFP_markov.3} and demonstrates Prop.
\ref{prop.3} for
any (classical) Markow process.

\section{Extension to fractional orders}\label{sec.fractional}

In the two previous sections (Sects.\ref{sect.sii_porcesses}-
\ref{sect.Markov}),
the fact that the indices $n \in J$ should be integers intervene at
best only in the correspondence between (integer order) differentiation
${\frac {\partial ^{n}}{\partial x^{n}}}$ (in Eq. \ref{eq_NFP_markov.3})
and powers of the conjugate variable $k^n$ (in Eq. \ref{eq.cumulant_gen}).
However, by the very definition
of fractional differentiation (e.g.\cite{Erdelyi}), this correspondence holds
also for non integer orders. However, there is not a unique definition of
fractional differentiation and therefore, as discussed in some details in
\cite{Yan}),
we cannot expect to have a unique expression of the Fractional Fokker-Planck
equation.

Since it will be sufficient for the following to consider
an expansion of the characteristic function involving fractional powers of
only the wave number $|k|$, it is interesting to  consider
Riesz's definition of a fractional differentiation.
Indeed, the latter corresponds to consider
fractional powers of the Laplacian:

\begin{equation}\label{eq_riesz}
-(-\Delta )^{\alpha /2}f(x)=F^{-1}[|k|^{\alpha}\hat{f}(k)]
\end{equation}

\noindent which has furthermore the advantage of being valid for the vector
cases.
However, we will see in Sect. \ref{sec.high.dim} that in general it does
not apply in a
straightforward manner for d-dimensional stable L\'evy motions. Indeed the
latter introduces
rather (one-dimensional) directional Laplacians, i.e. (one-dimensional)
Laplacians along a
given direction $\u{u}$ ($\mid \u{u} \mid =1$) :

\begin{equation} \label{eq_direct_laplacian}
-(-\Delta_{\u{u}})^{\alpha /2}f(x)=F^{-1}[|(\u{k},\u{u})|^{\alpha}\hat{f}(k)]
\end{equation}

\noindent where(.,.) denotes the scalar product. On the other hand, it will
be useful
to consider the fractional power of the contraction of the Laplacian tensor
$\u{\u{\Delta}}$:

\begin{equation} \label{eq_laplace_tensor}
\Delta_{i,j} = {\frac{\partial }{\partial x_i}} {\frac{\partial }{\partial
x_j}}
\end{equation}

\noindent by a tensor $\u{\u{\sigma }}$, with the following definition:

\begin{equation} \label{eq_sigma_laplacian}
-(-\u{\u{\Delta}}:\u{\u{\sigma}}.\u{\u{\sigma}}^{*})^{\alpha \over 2}
\equiv F^{-1} [\mid (\u{k}, \u{\u{\sigma }}.\u{\u{\sigma }}^{*}
.\u{k}\mid^{\alpha \over 2}]
=  F^{-1} [\mid \u{\u{\sigma }}^{*} .\u{k}\mid^{\alpha}]
\end{equation}

\section{Levy case}\label{Levy case}

The second characteristic function  of the increments $\delta L$
of the (scalar) Levy forcing is the following:

\begin{equation}  \label{eq_fc_dL}
\delta K_{L}(k,\delta t) =
\delta t [ik\gamma - D|k|^\a ( 1- i\beta {\frac{k}{|k|}} \omega(k, \alpha)
] + o(\delta t)
\end{equation}

\noindent where $\omega(k, \alpha)$ is defined by:

\begin{equation}  \label{eq_omega_k}
\alpha \neq1:~~ \omega (k,\alpha ) \equiv \omega (\alpha ) = tan {\frac{
\pi \alpha }{2}}; ~~~~\alpha = 1:~~\omega (k,\alpha ) ={\frac{\pi }{2}}
log|k|
\end{equation}

Considering an Ito-like forward integration of Eq.\ref{eq_langevin}, the
increments $\delta L$
generates the following (first) characteristic function for the increments
$\delta X$ of the motion
$X(t)$:

\begin{equation} \label{eq_fc_dZ}
\delta Z_{X}(k,\delta t|x- \delta x, t)  = e^{ik m(X,t)}
\delta Z_{\sigma L}(k,\delta t|x, t) + + o(\delta t)
\end{equation}

\noindent which yields the following elementary cumulant generating
function $\delta K_{X}$:

\begin{eqnarray}  \label{eq_fc_dK}
\delta K_{X}(k,\delta t|x,t)=
{\delta t}
 [ik m(x,t) + ik\gamma \sigma(x,t) \\
 - D|k|^\a (1- i\beta {\frac{k}{|k|}} \omega(k, \alpha)) \sigma(x,t)^{\alpha}]
+ o(\delta t)
\end{eqnarray}

\noindent and which is of the same type as Eq.\ref{eq.cumulant_gen}, with $J=
\{1,\alpha \}$. Therefore, as discussed in Sect.\ref{sec.fractional}, we
have fractional differentiations in the corresponding
Eq.\ref{eq_NFP_markov.3},
which will precisely correspond to Eq.\ref{eq_fokker_planck}, and therefore
establishes Prop. \ref{prop.1}.

Let us discuss briefly the regularity constraints that
should be satisfied by the nonlinear function $\sigma ((X(t),t)\geq 0$ and
$m((X(t),t)$.
Obviously, they should be measurable.
On the other hand, the uniqueness of the solution should require, as for the
classical nonlinear
Fokker-Planck equation (e;g. \cite{Oksendal}), a Lipschitz condition for both
$\sigma ((X(t),t)$ and $m((X(t),t)$.

\section{Extension to vector processes}\label{sec.high.dim}

With but one important exception, the extension of the previous
results to higher dimensions is rather straightforward. The starting
point of this extension is the following  nonlinear stochastic equation
($\u{X}(t) \in R^d$):

\begin{equation}  \label{eq_ddim_langevin}
d \u{X}(t) = \u{m}(\u{X}(t),t) dt + \u{\u{\sigma}}(\u{X}(t),t).d\u{L}
\end{equation}

\noindent where $\u{m}$ and $\u{\u{\sigma}}$
are the natural vector, respectively tensor, extensions
of the deterministic-like trend, respectively modulation of the random
driving force.
$\u{L}$ is a d-dimensional Levy stable motion and, as discussed below,
the expression of its characteristic function
corresponds to the source of the difficulty in extending the scalar results
to high dimensions.
On the contrary, it is straightforward to check that
Props. \ref{prop.2}, \ref{prop.3} are valid in the d-dimensional
case, with the following extensions ($\u{x} \in R^d$)  for Eq.
\ref{eq.FP.generation}:

\begin{equation} \label{eq.ddim.FP.generation}
{\frac{\partial p}{\partial t}}(\u{x},t|x_{0},t_{0}) =
\int {dy {\frac{\partial \widetilde{K}}{\partial t}}(\u{x}-\u{y}|\u{y},t)
p(\u{y},t|x_{0},t_{0})}
\end{equation}

\noindent and for Eq. \ref{eq_NFP_markov.3} ($\u{n}\in J \subseteq {\bf
N}^d, \mid \u{n} \mid = \sum_{i=1}^{d} n_i $):

\begin{equation} \label{eq_ddim_NFP_markov.3}
{\frac{\partial p}{\partial t}}(\u{x},t|\u{x}_{0},t_{0})=\sum_{\u{n} \in J}
  {\frac {\partial ^{\mid \u{n} \mid}}
  {\partial x_{1}^{n_1} \partial x_{2}^{n_2}..\partial x_{d}^{n_d}}}
[A_{\u{n}}(\u{x},t)  p(\u{x},t|\u{x}_{0},t_{0})]
\end{equation}

\noindent the relation to the cumulants $C_{\u{n}}$ of the increments is now:

\begin{equation}
A_{\u{n}}(\u{x},t) = {\frac{(-1)^{\mid \u{n} \mid}} {(n_1)!(n_2)!..(n_d)!}}
C_{\u{n}}(\u{x},t)
\end{equation}

On the other hand, Eq. \ref{eq_ddim_langevin} yields the following
extension to Eq.\ref{eq_fc_dZ}:

\begin{equation} \label{eq_ddim_fc_dZ}
\delta Z_{\u{X}}(\u{k},\delta t|\u{x}, t)  = e^{i~\u{k}.\u{m}(\u{x},t)}
\delta Z_{\u{\u{\sigma}}.\u{L}}(\u{k},\delta t|\u{x}, t)
\end{equation}

\noindent and therefore we have:

\begin{equation} \label{eq_ddim_fc_dK}
\delta K_{\u{X}}(\u{k},\delta t|\u{x}, t)  = i~\u{k}.\u{m}(\u{x},t) +
\delta K_{\u{L}}(\u{\u{\sigma}}^{*}.\u{k},\delta t|\u{x}, t) + o(\delta t)
\end{equation}

Let us recall that a stable L\'evy vector in the classical sense
\cite{Levy, Paulauskas, Nikias} (see \cite{Generalized_levy_vector}
for a discussion and a generalization) corresponds to the limit of a sum of
jumps,
with a power-law distribution, along
random directions $\u{u} \in \partial B_1$, $B_1$ being the unit ball,
distributed according to a (positive) measure $d\Sigma(\underline{u})$.
The latter, which generalizes the scale parameter $D$ of the scalar case,
is the source of the difficulty since in general the probability
distribution of a
stable L\'evy vector depends on this measure,
and therefore is a non parametric distribution. However, as discussed below,
there is at least a trivial exception: the case of isotropic stable L\'evy
vectors.

Corresponding to our previous remarks, a (classical) stable L\'evy vector
has the following (Fourier) cumulant generating function:

\begin{equation}\label{eq_ddim_classical_K}
K_{\underline{L}}(\underline{k})= \delta t [i
(\underline{k},\underline{\gamma})
-\int_{\underline{u} \in \partial B_1}
(i\underline{k},\underline{u})^{\alpha} d\Sigma(\underline{u}) ] + o (\delta t)
\end{equation}

\noindent which yields with the help of the Eq.\ref{eq_ddim_fc_dK}:

\begin{equation}\label{eq_ddim_classical_Ktilda}
\frac{\partial }{\partial t} {\widetilde{K}}_{\underline{X}}(\underline{k})=
-div (\underline{m} + \underline{\underline{\sigma}}.\underline{\gamma})
- F^{-1}[\int_{\underline{u} \in \partial B_1}
(i \underline{\underline{\sigma
}}^{*}(\underline{x},t).\underline{k},\underline{u})^{\alpha}
d\Sigma(\underline{u}) ]
\end{equation}

The scalar case (Eq.\ref{eq_fc_dL}) corresponds to:

\begin{equation} \label{eq_1dim_dSigma}
 0 \le p \le 1:\beta= 2p-1,  ~ d\Sigma(u) = D cos ({{\pi \alpha} \over
{2}})[p \delta_{(u-1)} + (1-p) \delta_{(u+1)}]
\end{equation}

For any dimension d, the second term on the right hand side of
Eq.\ref{eq_ddim_classical_Ktilda}
corresponds to a fractional differentiation
operator of order $\alpha$.  This operator can be slightly re-arranged.
With the help of
the odd $d\Sigma^{-}(\underline{u})$ and even $d\Sigma^{+}(\underline{u})$
parts of the measure $d\Sigma(\underline{u})$,

\begin{equation}
2~d\Sigma^{+}(\underline{u}) = d\Sigma(\underline{u}) +
d\Sigma(-\underline{u});~~~
2~d\Sigma^{-}(\underline{u}) = d\Sigma(\underline{u}) - d\Sigma(-\underline{u})
\end{equation}

\noindent and the identity ($\theta$ being the Heaviside function):

\begin{equation} \label{eq_identity}
(ik)^{\alpha} =  |k|^{\alpha }
[\theta (k) e^{i{\alpha \pi \over  2}} + \theta (-k) e^{-i{\alpha  \pi\over
2}}]
\end{equation}

\noindent one can write the extension of Eq.\ref{eq_fokker_planck} under
the following form:

\begin{eqnarray}  \label{eq_ddim_fokker_planck}
{\frac{\partial }{\partial t}} p(\underline{x},t|\underline{x}_0,t_0)
 =  - div[ \underline{m}(\underline{x},t)
+\underline{\underline{\sigma }}(\underline{x},t).\underline{\gamma}) ]
p(\underline{x},t|\underline{x}_0,t_0) ]  \nonumber \\
 -   [ <(-\u{\u{\Delta}}:\u{\u{\sigma}}.\u{\u{\sigma}}^{*})^{\alpha \over
2}>_{\Sigma^{+}}
-<(\u{\nabla}.\u{\u{\sigma}}^{*}).(-\u{\u{\Delta}}:\u{\u{\sigma}}.\u{\u{\sigma}}
^{*})^{{\alpha - 1} \over 2}>_{\Sigma^-}]
p(\underline{x},t|\underline{x}_0,t_0)
\end{eqnarray}

\noindent where the symmetric and antisymmetric operators are defined ,
similarly
to Eq.\ref{eq_sigma_laplacian}, in the following manner:

\begin{eqnarray} \label{eq_frac.operators}
- <(-\u{\u{\Delta}}:\u{\u{\sigma}}.\u{\u{\sigma}}^{*})^{\alpha \over
2}>_{\Sigma^{+}} &=&
\int_{\underline{u} \in \partial B_1} d\Sigma^+(\underline{u})
F^{-1} [ \mid (\u{\u{\sigma }}*(\u{x},t).\u{k},\u{u}) \mid^{\alpha}] \\
-
<(\u{\nabla}.\u{\u{\sigma}}^{*}).(-\u{\u{\Delta}}:\u{\u{\sigma}}.\u{\u{\sigma}}^
{*})^{{\alpha - 1} \over 2}>_{\Sigma^-} &=&
\int_{\underline{u} \in \partial B_1} d\Sigma^-(\underline{u})
F^{-1} [ (- i \u{\u{\sigma }}^{*}(\u{x},t).\u{k},\u{u})\mid (\u{\u{\sigma
}}^{*}(\u{x},t).\u{k},\u{u})
\mid^{\alpha -1}]
\end{eqnarray}

In general, each operator corresponds to a rather complex integration
(which is indicated
by the symbol $<.>_{\Sigma}$)  of directional fractional Laplacians
(Eq.\ref{eq_direct_laplacian}).
However, the symmetric operator becomes simpler as soon as the even part
$d\Sigma^+$  of the measure $d\Sigma$
is isotropic. Indeed, the integration over directions yields only a
prefactor $D$:

\begin{eqnarray} \label{eq_sym.frac.operator}
<-(\u{\u{\Delta}}:\u{\u{\sigma}}.\u{\u{\sigma}}^{*})^{\alpha \over
2}>_{\Sigma^{+}} =
D~ (-\u{\u{\Delta}}:\u{\u{\sigma}}.\u{\u{\sigma}}^{*})^{\alpha \over 2}
\nonumber \\
D = \int_{\underline{u} \in \partial B_1} d\Sigma^{+}(\underline{u}) \mid
(\u{u}_{1},\u{u}) \mid^{\alpha}
\end{eqnarray}

\noindent and for $\alpha =2$ this corresponds to the classical term
($\u{\u{\Delta}}:\u{\u{\sigma}}.\u{\u{\sigma}}^{*}$)
of the standard d-dimensional Fokker-Planck equation.  If $d\Sigma$ itself
is rotation invariant,
then the asymmetric operator vanishes, since $d\Sigma^- = 0$. If
furthermore, $\u{\u{\sigma}}$
is rotation invariant, i.e. $\u{\u{\sigma}}= \sigma \u{\u{1}}$,
then one obtains the following Fractional Fokker-Planck equation:

\begin{eqnarray}  \label{eq_ddim_sym_fokker_planck}
{\frac{\partial }{\partial t}} p(\underline{x},t|\underline{x}_0,t_0)
&=& - div[ \underline{\underline{\sigma }} . \underline{\gamma}
(\underline{x},t) + \underline{m}(\underline{x},t)]
p(\underline{x},t|\underline{x}_0,t_0)  \\
&-& D~[(-\Delta )^{\alpha /2}] \sigma (x,t)^{\alpha}
p(\underline{x},t|\underline{x}_0,t_0)
\end{eqnarray}

Therefore, as one might expect it, due to the rotation symmetries,
this corresponds to a rather trivial extension
of the standard gaussian case: a fractional power $\alpha$ of the
d-dimensional Laplacian,
as in the pure scalar case (Eq.\ref{eq_fokker_planck}).  Obviously, the
integration performed
in Eq.\ref{eq_ddim_fokker_planck} is also greatly
simplified as soon as $d\Sigma(\underline{u})$ is discrete, i.e. its support
corresponds to a discrete set of directions $\u{u}_i$.

On the other hand, let us note that the framework of generalized stable
L\'evy vectors \cite{Generalized_levy_vector}, allows one to introduce
a much stronger anisotropy than the the measure $d\Sigma$ does if for
classical stable L\'evy vectors. This therefore diminishes the importance of
the asymmetry of the latter.
Indeed, the components of a generalized stable L\'evy vector do not
have necessarily the same L\'evy stability index, the latter being
generalized into
a second rank tensor. Similarly, the differential operators involved in the
corresponding
Fractional Fokker-Planck equation  have no longer a unique order of
differentiation.
This is rather easy to check in case of a discrete measure
$d\Sigma(\underline{u})$
and we will explore elsewhere the general case.

\section{Conclusion}\label{Conclusion}

We have derived a Fractional Fokker-Planck equation, i.e. a kinetic equation
which involves fractional derivatives, for the evolution of the
probability distribution of nonlinear stochastic differential equations
driven by non-Gaussian Levy stable noises. We first established this equation
in the scalar case, where it has a rather compact expression
with the help of fractional powers of the Laplacian, and then discussed
its extension to the vector case.
This  Fractional Fokker-Planck equation generalizes
broadly previous results obtained for a linear Langevin-like equation with
a L\'evy forcing, as well as the standard Fokker-Planck equation for
a nonlinear Langevin equation with a Gaussian forcing.

\section{Acknowledgments}\label{Acknowledgments}

We would like to thank Dr. James Brannan for helpful
discussions. Part of this work was performed while
Daniel Schertzer was visiting Clemson University.


\end{document}